\documentclass[12pt,dvips]{article}

\usepackage{amsmath,amsfonts,amssymb,amsthm,graphicx,verbatim,subfig}
\usepackage[all]{xy}

\ifx\pdftexversion\undefined

\usepackage[a4paper,colorlinks,link
=black,filecolor=black,citecolor=black,urlcolor=black,pdfstartview=FitH]{hyperref}
\else

\usepackage[a4paper,colorlinks,linkcolor=black,filecolor=black,citecolor=black,urlcolor=black,pdfstartview=FitH]{hyperref}
\fi

\font\sixbb=msbm6
\font\eightbb=msbm8
\font\twelvebb=msbm10 scaled 1095
\newfam\bbfam
\textfont\bbfam=\twelvebb \scriptfont\bbfam=\eightbb
                           \scriptscriptfont\bbfam=\sixbb

\newtheorem{theorem}{\bf Theorem}[section]
\newtheorem{claim}[theorem]{\bf Claim}

\newtheorem{proposition}[theorem]{\bf Proposition}
\newtheorem{corollary}[theorem]{\bf Corollary}

\newcommand{\enp}{\begin{flushright} $\Box$ \end{flushright}}
\newcommand{\beq}[0]{\begin{equation}}
\newcommand{\enq}[0]{\end{equation}}
\newcommand{\dn}{\Delta_{n-1}}

\newcommand{\prob}{{\rm Pr}}

\newcommand{\cg}{{\cal G}}

\newcommand{\namedref}[2]{\hyperref[#2]{#1~\ref*{#2}}}

\title{Bounded Quotients of the Fundamental Group of a Random $2$-Complex}
\begin{document}
\author{Roy Meshulam\thanks{Department of Mathematics,
Technion, Haifa 32000, Israel. e-mail:
meshulam@math.technion.ac.il~. Supported by an ISF grant. }
}

\maketitle
\pagestyle{plain}
\begin{abstract}

Let $\dn$ denote the $(n-1)$-dimensional simplex. Let $Y$ be a
random $2$-dimensional subcomplex of $\dn$ obtained by starting
with the full $1$-skeleton of $\dn$ and then adding
each $2$-simplex independently with probability $p$. For a fixed $c>0$ it is shown that if
$p=\frac{(6+7 c)\log n}{n}$ then a.a.s. the fundamental group $\pi_1(Y)$
does not have a nontrivial quotient of order at most $n^c$.

\end{abstract}

\section{Introduction}
\ \ \
For a simplicial complex $X$, let $X^{(k)}$ denote the $k$-skeleton of $X$ and let
$f_k(X)$ be the number of $k$-simplices of $X$.
Let $\dn$ be the $(n-1)$-dimensional simplex on the vertex set
$V=[n]$. Let $Y(n,p)$ denote the probability space of complexes
$\dn^{(1)} \subset Y \subset \dn^{(2)}$ with
probability measure
$$\Pr(Y)=p^{f_2(Y)}(1-p)^{\binom{n}{3}-f_2(Y)}~.$$

The threshold probability for the vanishing of the first homology with fixed finite abelian coefficient group was determined in \cite{LM06,MW09}. Let $\omega(n)$ be an arbitrary function that tends to infinity with $n$.
\begin{theorem}[\cite{LM06,MW09}]
\label{gen}  Let $R$ be a fixed finite abelian group. Then
$$
\lim_{n \rightarrow \infty} \prob ~[~Y \in Y(n,p):
H_1(Y;R)=0 ~]= \left\{
\begin{array}{ll}
        0 & p=\frac{2\log n -\omega(n)}{n} \\
        1 & p=\frac{2\log n+ \omega(n)}{n}~.
\end{array}
\right.
$$
\end{theorem}

The  threshold probability for the vanishing of the fundamental group was determined by
Babson, Hoffman and Kahle \cite{BHK}.
\begin{theorem}[\cite{BHK}]
\label{bhkt}
Let $\epsilon>0$ be fixed, then
$$
\lim_{n \rightarrow \infty} \prob ~[~Y \in Y(n,p):
\pi_1(Y)=0 ~]= \left\{
\begin{array}{ll}
        0 & p=n^{-\frac{1}{2}-\epsilon} \\
        1 & p=(\frac{3 \log n +\omega(n)}{n})^{1/2}~.
\end{array}
\right.
$$
\end{theorem}
In view of the gap between the thresholds for the vanishing of $H^1(Y;R)$ ($R$ finite) and
for the triviality of $\pi_1(Y)$, Eric Babson (see problem (8) on page 58 in \cite{FGS})
asked what is the threshold probability such that a.a.s. $\pi_1(Y)$ does not have a quotient equal to some finite group. Addressing Babson's question we prove the following
\begin{theorem}
\label{nonab}
Let $c >0$ be fixed and let $p=\frac{(6+7 c)\log n}{n}$. Then a.a.s. $\pi_1(Y)$
does not contain a nontrivial normal subgroup of index at most $n^c$.
\end{theorem}
\noindent
{\bf Remark:} The constant $6+7c$ may be improved using a more careful analysis as in \cite{LM06,MW09}. For example, for any {\it fixed} non-trivial finite group  $G$, if
$p=\frac{2 \log n+\omega(n)}{n}$ then a.a.s. $G$ is not a homomorphic image of $\pi_1(Y)$.

The proof of Theorem \ref{nonab} is an adaptation of an argument in \cite{LM06,MW09} to the non-abelian setting.
In Section \ref{naco} we recall the notion of non-abelian first cohomology and its relation with the fundamental group.
In Section \ref{expsim} we compute the expansion of the $(n-1)$-simplex. The results of Sections \ref{naco} and \ref{expsim} are used
in section \ref{fmw} to prove Theorem \ref{nonab}.

\section{Non-abelian first cohomology}
\label{naco}
\ \ \
Let $X$ be a simplicial complex and let $G$ be a multiplicative group.
We recall the definition of the first cohomology $H^1(X;G)$ of $X$ with $G$ coefficients (see e.g. \cite{Olum58}).
For $0 \leq k \leq 2$ let $X(k)$ be the set of all ordered $k$-simplices of $X$.
Let $C^0(X;G)$ denote the group of  $G$-valued functions on $X(0)$ with pointwise multiplication,
and let $$C^1(X;G)=\{\phi:X(1) \rightarrow G: \phi(u,v)=\phi(v,u)^{-1}\}.$$
The $0$-coboundary operator $d_0:C^0(X;G) \rightarrow C^1(X;G)$ be given by
$$d_0 \psi(u,v)=\psi(u)\psi(v)^{-1}.$$
For $\phi \in C^1(X;G)$ and $(u,v,w) \in X(2)$ let
$$d_1 \phi (u,v,w)=\phi(u,v)\phi(v,w)\phi(w,u).$$
The set of $G$-valued $1$-cocycles of $X$ is given by
$$Z^1(X;G)=\{\phi \in C^1(X;G):d_1\phi(u,v,w)=1 {\rm ~for~all~}(u,v,w) \in X(2)\}.$$
Define an action of $C^0(X;G)$ on $C^1(X;G)$ as follows.
For $\psi \in C^0(X;G)$ and  $\phi \in C^1(X;G)$ let
$$\psi . \phi (u,v)=\psi(u) \phi(u,v) \psi(v)^{-1}.$$
Note that $d_0 \psi =\psi . 1$ and that $Z^1(X;G)$ is invariant under the action of $C^0(X;G)$.
For $\phi \in C^1(X;G)$ let $[\phi]$ denote the orbit of $\phi$ under the action of $C^0(X;G)$.
The first cohomology of $X$ with coefficients in $G$ is the set of orbits
$$H^1(X;G)=\{[\phi]: \phi \in Z^1(X;G) \}.$$

Let $Hom(\pi_1(X),G)$ denote the set of homomorphisms from $\pi_1(X)$ to $G$.
Let $G$ act on $Hom(\pi_1(X),G)$ by conjugation and for $\varphi \in Hom(\pi_1(X),G)$ let $[\varphi]$ denote the orbit of $\varphi$ under this action.
Let $$Hom(\pi_1(X),G)/G=\{[\varphi]: \varphi \in Hom(\pi_1(X),G)\}.$$
The following observation is well known (see (1.3) in \cite{Olum58}). For completeness we outline a proof.
\begin{claim}
\label{fhpi}
For $\dn^{(1)} \subset X \subset \dn^{(2)}$
there is a bijection $$\mu:Hom(\pi_1(X),G)/G \rightarrow H^1(X;G)$$ that maps $[1] \in Hom(\pi_1(X),G)/G$
to $[1] \in H^1(X;G)$.
\end{claim}
\noindent
{\bf Proof:}
We identify $\pi_1(X)$ with the group generated by $\{e_{ij}: 2 \leq i \neq j \leq n\}$
modulo the relations
\begin{itemize}
\item
$e_{ij}e_{ji}=1$.
\item
$e_{ij}=1$ if $(1,i,j) \in X(2)$.
\item
$e_{ij}e_{jk}e_{ki}=1$ if $(i,j,k) \in X(2)$.
\end{itemize}
For $\varphi \in Hom(\pi_1(X),G)$ let $F(\varphi) \in C^1(X;G)$ be given by
$$
F(\varphi)(i,j)=
\left\{
\begin{array}{ll}
        \varphi(e_{ij}) & 2 \leq i \neq j \leq n \\
        1 & otherwise.
\end{array}
\right.
$$
It can be checked that $F(\varphi) \in Z^1(X;G)$ and that the mapping $$\tilde{F}:Hom(\pi_1(X),G)/G \rightarrow H^1(X;G)$$
given by $\tilde{F}([\varphi])=[F(\varphi)]$ is the required bijection.
{\enp}
\noindent
In particular we obtain the following
\begin{corollary}
\label{hhpp}
Let $N>1$. Then
$\pi_1(X)$ contains a nontrivial normal subgroup of index at most $N$  iff  $H^1(X;G) \neq \{[1]\}$ for some nontrivial simple group $G$ of order  at most $N$.
{\enp}
\end{corollary}

\section{$1$-Expansion of the Simplex}
\label{expsim}
Let $\phi \in C^1(\dn;G)$.
The {\it support size} of $\phi$ is
$$\|\phi\|=|\{u,v\} \in \binom{[n]}{2}: \phi(u,v) \neq 1\}|.$$
The {\it weight} of the orbit $[\phi]$ is
$$\|[\phi]\|=\min\{\|\psi . \phi\|: \psi \in C^0(\dn;G)\}.$$
Let $B(\phi)$ be the support of $d_1\phi$, i.e.
$$
B(\phi)=\left\{ \{u,v,w\} \in \binom{[n]}{3} : d_1\phi(u,v,w)\neq 1 \right\}
$$
and let $\|d_1\phi\|=|B(\phi)|$.
The following result is an adaptation of Proposition 3.1 of \cite{MW09} to the non-abelian setting.
\begin{proposition}
\label{bw1}
Let $\phi \in C^1(\dn;G)$ then
$$\|d_1\phi\| \geq \frac{n \|[\phi]\|}{3}.$$
\end{proposition}
\noindent
{\bf Proof:} For $u \in \dn(0)$ define $\phi_u \in C^0(\dn;G)$ by
$$\phi_u(v)=
\left\{
\begin{array}{ll}
        1 & v=u \\
        \phi(u,v) & v \neq u .
\end{array}
\right.
$$
Note that if $(u,v,w) \in \dn(2)$ then
\begin{equation*}
\begin{split}
d_1\phi(u,v,w)&=\phi(u,v)\phi(v,w)\phi(w,u) \\
&=\phi_u(v)\phi(v,w)\phi_u(w)^{-1}=\phi_u.\phi(v,w).
\end{split}
\end{equation*}
Therefore
\begin{equation*}
\begin{split}
6\|d_1\phi\|&=|\{(u,v,w) \in \dn(2):d_1\phi(u,v,w) \neq 1\}| \\
&= |\{(u,v,w) \in \dn(2): \phi_u.\phi(v,w) \neq 1\}| \\
&=\sum_{u} 2\|\phi_u . \phi\| \geq 2n \|[\phi]\|.
\end{split}
\end{equation*}
{\enp}

\section{Proof of Theorem \ref{nonab}}
\label{fmw}

Let $G$ be a finite group. For a subcomplex $\dn^{(1)} \subset X \subset \dn^{(2)}$ we identify $H^1(X;G)$ with its image under the natural injection
$H^1(X;G) \rightarrow H^1(\dn^{(1)};G)$. If $\phi \in C^1(\dn;G)$ then
$[\phi] \in H^1(X;G)$ iff $d_1\phi(u,v,w)=1$ whenever $(u,v,w) \in X(2)$.
It follows that in the probability space $Y(n,p)$
$$
\prob\left[ [\phi] \in H^1(Y;G)\right] =(1-p)^{\|d_1\phi\|}.$$
Therefore
\begin{equation}
\label{unbod}
\begin{split}
\prob \left[H^1(Y;G) \neq \{[1]\}\right] &\leq \sum_{[1] \neq [\phi] \in  H^1(\dn^{(1)};G)} \prob\left[
[\phi] \in H^1(Y;G)\right] \\
&= \sum_{[1] \neq [\phi] \in  H^1(\dn^{(1)};G)}(1-p)^{\|d_1 \phi\|}~.
\end{split}
\end{equation}
Suppose now that is $|G| \leq n^c$. Then
by (\ref{unbod}) and Proposition \ref{bw1}
\begin{equation}
\label{unbd1}
\begin{split}
& \prob \left[H^1(Y;G) \neq \{[1]\}\right] \\ &\leq
\sum_{k \geq 1} \sum_{\|[\phi]\|=k} (1-p)^{\frac{kn}{3}} \\
&\leq \sum_{k \geq 1} \binom{\binom{n}{2}}{k} |G|^k  (1-\frac{(6+7c)\log n}{n})^{\frac{kn}{3}} \\
&\leq \sum_{k \geq 1} n^{2k} n^{ck} n^{-\frac{(6+7c)k}{3}}=O(n^{-\frac{4c}{3}}).
\end{split}
\end{equation}

Let $\cg(N)$ be the set of all nontrivial simple groups of size at most $N$.
As there are at most $2$ non-isomorphic simple groups of the same order, it follows that $|\cg(N)| \leq 2N$.
Combining Corollary \ref{hhpp} and (\ref{unbd1}) it follows that the probability that $\pi_1(Y)$ contains a nontrivial normal subgroup of index $ \leq n^c$ is at most
$$\sum_{G \in \cg(n^c)} \prob[Y \in Y(n,p)~:~H^1(Y;G) \neq \{[1]\}] \leq $$
$$O(|\cg(n^c)|n^{-\frac{4c}{3}})=O(n^{-\frac{c}{3}}).$$
{\enp}

\ \\ \\
{\bf Acknowledgement:} This work was carried out while the author was attending Oberwolfach Arbeitsgemeinschaft on Topological Robotics
on October 2010. The author would like to thank Eric Babson for raising the question and for his helpful comments.

\end{document}